\documentclass[12pt]{article}

\usepackage{amsmath}
\usepackage{amsfonts}
\usepackage{amssymb}
\usepackage{graphicx, color}
\usepackage{hyperref}
\usepackage{epsfig}
\usepackage{epstopdf} 

\definecolor{darkgreen}{rgb}{0,0.55,0}

\newtheorem{proposition}{Proposition}[section]
\newtheorem{theorem}{Theorem}[section]
\newtheorem{lemma}[theorem]{Lemma}
\newtheorem{corollary}[theorem]{Corollary}

\DeclareSymbolFont{AMSb}{U}{msb}{m}{n}
\DeclareMathSymbol{\N}{\mathbin}{AMSb}{"4E}
\DeclareMathSymbol{\Z}{\mathbin}{AMSb}{"5A}
\DeclareMathSymbol{\R}{\mathbin}{AMSb}{"52}
\DeclareMathSymbol{\Q}{\mathbin}{AMSb}{"51}
\DeclareMathSymbol{\I}{\mathbin}{AMSb}{"49}

\begin{document}

\title{Uniqueness of solutions of mean field equations in $\R^2$ }
\author{{Changfeng Gui\footnote{Department of Mathematics, University of Texas at San Antonio, Texas, USA. E-mail: changfeng.gui@utsa.edu.  }
\qquad Amir Moradifam\footnote{Department of Mathematics, University of California, Riverside, California, USA. E-mail: moradifam@math.ucr.edu. }}}
\date{\today}

\smallbreak \maketitle

\begin{abstract}
In this paper, we prove uniqueness of  solutions of mean field equations with  general boundary conditions for the critical and subcritical 
total mass regime, extending the earlier results for  null Dirichlet boundary condition.   The proof is based on new Bol's inequalities for weak radial solutions obtained from rearrangement  of the solutions.

\end{abstract}
\maketitle

\section{Introduction} 
Let $\Omega \subset \R^2$ be an open bounded domain and  consider the mean field equation 

\begin{eqnarray}\label{MainEquation0}
\left\{ \begin{array}{ll}
\Delta u+\rho \frac{e^u}{\int_{\Omega}e^{u}}=0 &\text{in } \Omega\\
u=0 &\text{on }\partial \Omega.
\end{array} \right.
\end{eqnarray}
Suzuki \cite{Suzuki} proved that if $\Omega$ is simply-connected, then for $0<\rho<8\pi$  the equation (\ref{MainEquation0})  has a unique solution. Later in \cite{CCL} the authors extended this result to the case $\rho=8\pi$. Recently in \cite{BL} Bartolucci and Lin  extended the result to multiply connected domains. Indeed they proved the following. \\ \\
{\bf Theorem A.} (Theorem 2 in \cite{BL})  \textit{Let $\Omega$ be an open, bounded, and multiply connected domain of class $C^1$. Then  equation \eqref{MainEquation0} admits at most one solution for $0<\rho \leq 8\pi.$} \\  \\
The proof relies on a generalization of the classical  Bol's inequality for multiply-connected domains (see Theorem C below). A necessary and sufficient condition for the existence of a solution at the critical parameter $\rho=8\pi$ is also provided in \cite{BL}. 

In this paper, among other results, we study uniqueness of solutions of the general mean field equation 
\begin{eqnarray}\label{MainEquationGeneral}
\left\{ \begin{array}{ll}
\Delta u+\rho \frac{K(x)e^u}{\int_{\Omega}K(x)e^{u}}=f &\text{in } \Omega\\
u=g &\text{on }\partial \Omega,
\end{array} \right.
\end{eqnarray}
on simply connected domains, where $K$ is a prescribed positive $C^2$  function. We shall prove the following uniqueness results for $\rho \le 8 \pi$.

\begin{theorem}\label{MainResult}
Let $\Omega$ be an open, bounded and simply-connected domain, and let $K \in C^2(\Omega)\cap  C(\overline{\Omega})$ be positive.  Assume that $v_i \in C^2(\Omega)\cap  C(\overline{\Omega})$, $i=1,2$, satisfy
\begin{equation}
\Delta v_i + Ke^{v_i}=f_i(x),
\end{equation}
where  $f_2 \geq f_1\geq -\Delta \ln(K)$ in $\Omega$. If $v_1 \not \equiv v_2$, $v_2-v_1=c$ on $\partial \Omega$ for some $c \in \R$, and 
\[\int_{\Omega} Ke^{v_1}dx=\int_{\Omega}K e^{v_2}dx=\rho,\]
then $\rho > 8\pi $. 
\end{theorem}

The above theorem is equivalent to the next uniqueness result. Indeed Theorem \ref{MainResult} follows from Theorem \ref{New-8Pi-LowerBound} by letting $w_i=\ln K+v_i$, i=1,2, and Theorem \ref{New-8Pi-LowerBound} follows from Theorem \ref{MainResult} by letting $K\equiv 1$. 

\begin{theorem}\label{New-8Pi-LowerBound}
Let $\Omega$ be an open, bounded, and simply-connected domain.  Assume that $w_i \in C^2(\Omega)\cap  C(\overline{\Omega})$, $i=1,2$, satisfy
\begin{equation}
\Delta w_i + e^{w_i}=f_i(x),
\end{equation}
where  $f_2 \geq f_1\geq 0$ in $\Omega$. If $w_1 \not \equiv w_2$, $w_2-w_1=c$ on $\partial \Omega$ for some $c \in \R$, and 
\[\int_{\Omega} e^{w_1}dx=\int_{\Omega}e^{w_2}dx=\rho,\]
then $\rho > 8\pi $. 

\end{theorem}
\vspace{.2cm}

\begin{corollary}\label{meainFieldUniqueness}
Let $\Omega$ be an open, bounded and simply-connected domain, and let $K \in C^2(\Omega)\cap  C(\overline{\Omega})$ be positive.   Assume that $u_i \in C^2(\Omega)\cap  C(\overline{\Omega})$, $i=1,2$ satisfy
\begin{equation}
\Delta u_i +\rho \frac{Ke^{u_i}}{\int_{\Omega}Ke^{u_i}}=f_i(x),
\end{equation}
where $f_2\geq f_1\geq -\Delta \ln (K)$ in $\Omega$, and $u_1\not \equiv u_2$. If $u_2-u_1=c$ on $\partial \Omega$ for some $c \in \R$, then $\rho > 8\pi $. 

\end{corollary}
\vspace{.3cm}

We also present the following uniqueness results on multiply-connected domains.

\begin{theorem}\label{New-8Pi-LowerBoundM}
Let $\Omega$ be an open, bounded and multiply-connected domain, and let $K \in C^2(\Omega)\cap  C(\overline{\Omega})$ be positive.   Assume that $v_i \in C^2(\Omega)\cap  C(\overline{\Omega})$, $i=1,2$, satisfy
\begin{equation}
\Delta v_i + Ke^{v_i}=f_i(x),
\end{equation}
where $f_2 \geq f_1\geq -\Delta \ln (K)$. If $v_1 \not \equiv v_2$, $v_1+\ln(K)=c_1$ and $v_2+\ln(K)=c_2$ on $\partial \Omega$ for some $c_1, c_2 \in \R$, and 
\[\int_{\Omega} Ke^{v_1}dx=\int_{\Omega}K e^{v_2}dx=\rho,\]
then $\rho > 8\pi $.
\end{theorem}

\begin{corollary}\label{meainFieldUniquenessM}
Let $\Omega$ be an open, bounded and multiply-connected domain, and let $K \in C^2(\Omega)\cap  C(\overline{\Omega})$ be positive.     Assume that $u_i \in C^2(\Omega)\cap  C(\overline{\Omega})$, $i=1,2$ satisfy
\begin{equation}
\Delta u_i +\rho \frac{K e^{u_i}}{\int_{\Omega}K e^{u_i}}=f_i(x),
\end{equation}
where where $f_2 \geq f_1\geq -\Delta \ln (K)$. If $u_1 \not \equiv u_2$, $u_1+\ln(K)=c_1$ and $u_2+\ln(K)=c_2$ on $\partial \Omega$ for some $c_1, c_2 \in \R$, then $\rho > 8\pi $. 
\end{corollary}

We should mention that Theorems 1.1 and 1.2, and Corollary 1.3 are known in some special cases and also when the the wight $K$ is singular (see \cite{BL0}).

We shall use  a new idea  from \cite{G-M1} regarding the  rearrangement  of the solutions according to the standard metric on a sphere  (projected to $\R^2)$  and  compare the total masses of the solutions.  In particular,  we  shall show  a reversed  Bol's  inequality  in exterior domain for weak radial solutions (Proposition \ref{RadialBolPropExterior}).  

\section{Preliminaries}

Bol's isoperimetric inequality plays a crucial role in the proof of our main results. In this section we gather some results on Bol's inequality that will be used in subsequent sections. Let us first recall the classical Bol's isoperimetric inequality, see \cite{Bandle, BL0, Bol, CCL, LL, Suzuki}, and \cite{BZ} for a detailed history of the Bol's inequality. \\ \\
{\bf Theorem B.}
\textit{Let $\Omega \subseteq \R^2$ be a simply-connected and assume  $u \in C^2(\Omega)\cap C(\overline{\Omega})$ satisfies 
\begin{equation}\label{condition}
\Delta u+e^u\geq 0, \ \ \ \int_{\Omega}e^u dx \leq 8 \pi.
\end{equation}
Then for every $\omega \Subset \Omega$ of class $C^1$ the following inequality holds
\begin{equation}\label{BolInequality0}
\left( \int_{\partial \omega}e^{\frac{u}{2}}\right)^2\geq \frac{1}{2}\left( \int_{\omega} e^u\right) \left( 8 \pi-\int_{\omega}e^u\right).
\end{equation}
Moreover   the inequality in (\ref{BolInequality0}) is strict if $\Delta u+e^{u}>0$ somewhere in $\omega$ or $\omega$ is not simply-connected.  \\}

For $\lambda>0$ the function $U_\lambda$ defined by 
\begin{equation}\label{ULambda}
U_{\lambda}:=-2\ln(1+\frac{\lambda^2 |y|^2 }{8})+2\ln(\lambda)
\end{equation}
satisfies
\begin{equation*}
\Delta U_\lambda+e^{U_\lambda}=0,
\end{equation*}
and 
\[\left( \int_{\partial B_r}e^{\frac{ U_\lambda}{2}}\right)^2= \frac{1}{2}\left( \int_{B_r} e^{ U_\lambda}\right) \left( 8 \pi-\int_{B_r}e^{ U_\lambda}\right),\]
for all $r>0$ and $\lambda>0$, where $B_r$ denotes the ball of radius $r$ centered at the origin in $\R^2$. 

Note that Theorem B requires $\Omega$ to be simply-connected but $\omega$ can be multiply-connected. Recently in \cite{BL} Theorem B 
is extended to the case where  $\Omega \subset \R^2$  is  multiply-connected and  $u$ is constant on $\partial \Omega$.  \\ \\
{\bf Theorem C.} (Theorem 3 in \cite{BL})
\textit{Let $\Omega$ be an open and bounded domain of class $C^1$ in $\R^2$ and assume $u \in C^2(\Omega)\cap C(\overline{\Omega})$ satisfies  \eqref{condition}  and $u=c$ on $\partial \Omega$, for some constant $c\in \R$.  Then \eqref{BolInequality0} holds
for every $\omega \Subset \Omega$.  Moreover the inequality is strict if  $\Delta u+e^{u}>0$ somewhere in $\omega$ or $\omega$ is not simply-connected. } \\ 

Let $\Omega$ be an open, bounded, and multiply-connected domain of class $C^1$ in $\R^2$, and $\overline{\Omega^*}$  be the closure of the union of the bounded components of $\R^2 \setminus \partial \Omega$ and $\Omega^* = \overline{\Omega^*}\setminus  \partial \overline{\Omega^*}$. It is easy to see that $\Omega\subset \Omega^*$. Suppose $g\in C(\partial \Omega)$ satisfies
\begin{equation}\label{g}
g=G|_{\partial \Omega},
\end{equation}
where $G$ is Lipschitz continuous in $\overline{\Omega^*}$, $G$ is subharmonic in $\overline{\Omega^*}$ and harmonic in $\Omega$. The following more general result is also proved in \cite{BL}. \\ \\
{\bf Theorem D.} (Theorem 4 in \cite{BL})
\textit{Let $\Omega$ be an open, bounded, and multiply-connected domain of class $C^1$ in $\R^2$. Suppose $u \in C^2(\Omega)\cap C(\overline{\Omega})$ satisfies \eqref{condition} with $u=g$ on $\partial \Omega$, and $g\in C(\partial \Omega)$ satisfies (\ref{g}). Then \eqref{BolInequality0} holds 
for every $\omega \Subset \Omega$.  Moreover the inequality is strict  if   $\Delta u+e^{u}>0$ somewhere in $\omega$  or  $\omega$ is not simply-connected. } \\ \\

Next we shall recall some facts about rearrangements according to the metric on $\R^2$ which is the  stereographic projection of  the standard metric on the unit sphere. Such rearrangments are discussed in detail in \cite{G-M1}, but we also include it here for the sake of the readers. Let $\Omega \subset \R^2$ and $\lambda>0$, and suppose that $u\in C^2(\overline{\Omega})$ satisfies
\[\Delta u+e^{u}\geq 0.\]
Then any function $ \phi \in C^2(\overline{\Omega})$ which is constant on $\partial \Omega$ can be equimeasurably rearranged with respect to the measures $e^udy$ and $e^{U_\lambda}dy$ (see \cite{Bandle}, \cite{BL0}, \cite{CCL},  \cite{LL}, \cite{Suzuki}), where $U_\lambda$ is defined in (\ref{ULambda}). More precisely, for $t>\min_{y\in \overline\Omega}\phi$ define
\[\Omega_t:=\{\phi>t\} \subset \Omega,\]
and define $\Omega^*_t$ be the ball centered at the  origin in $\R^2$ such that 
\[\int_{\Omega^*_t}e^{U_{\lambda}}dy=\int_{\Omega_t}e^{u}dy:=a(t).\]
Then $a(t) $ is  a right-continuous function,  and $\phi^*: \Omega^* \rightarrow \R$ defined by $\phi^*(y):=\sup \{t\in \R: y\in \Omega^*_t\}$ provides an equimeasurable rearrangement of $\phi$ with respect to the measure $e^{u}dy$ and $e^{U_{\lambda}}dy$, i.e. 
\begin{equation}\label{rearrang}
\int_{\{\phi^*>t\}}e^{U_{\lambda}}dy=\int_{\{\phi>t\}}e^{u}dy, \ \ \forall  t>\min_{y\in \overline \Omega}\phi.
\end{equation}
We shall need the following lemma. 

\begin{proposition}\label{LipCont}
Let $u, \varphi \in C^1(\bar{\Omega})$ and assume that $\phi$ is constant on $\partial \Omega$. Let $\phi^*(r)$ be the equimeasurable rearrangement of $\phi$ with respect to the measure $e^{u}dy$ and $e^{U_{\lambda}}dy$. Then $\phi^*$ is Lipschitz continuous on $(\epsilon,R-\epsilon)$, for every $\epsilon>0$, where $R$ is the radius of $\Omega^*$.\\ 
\end{proposition}
{\bf Proof.} First note that the function $\phi^*$ is decreasing and the set 
\[\mathcal{T}:=\{t \geq \min_{\bar{\Omega}}\phi: \ \ (\phi^*)^{-1}(t) \ \ \hbox{is not a singleton}\}\]
has Lebesgue measure zero. Indeed $ (\phi^*)^{-1}(t)$ is a connected closed interval for all $t\in \mathcal{T}$. Let $0< r_1<r_2 < R$ and
\[a(t)=\int_{\{\phi^*>t\}}e^{U_{\lambda}}dy=\int_{\{\phi>t\}}e^{u}dy, \ \ \forall  t>\min_{y\in \overline \Omega}\phi.\]
For $\phi^*(r_1),\phi^*(r_2)\not \in \mathcal{T}$, we have 
\begin{eqnarray*}
a(\phi^*(r_2))-a(\phi^*(r_1))&=&\int_{\{\phi^*(|y|)>\phi^*(r_2)\}}e^{U_{\lambda}}dy-\int_{\{\phi^*(|y|)>\phi^*(r_2)\}}e^{U_{\lambda}}dy\\
&=&\int_{\{\phi(y)>\phi^*(r_2)\}}e^{u}dy-\int_{\{\phi(y)>\phi^*(r_2)\}}e^{u}dy\\
&=&\int_{\{\phi^*(r_2)\leq \phi(y)\leq \phi^*(r_1)\}}e^{u}dy \\
&=&\int_{\{\phi^*(r_2)\leq \phi^*(|y|)\leq \phi^*(r_1)\}}e^{U_\lambda}dy. 
\end{eqnarray*}
Now let $m:=\min \limits_{\overline{\Omega}} e^{u(y)}$, $M_1:=\max \limits_{\overline{\Omega^*}}e^{U_\lambda(y)}$, and $M_2:=\max \limits_{\overline{\Omega}}|\nabla \phi|$. Then it follows from the above equality that  
\begin{eqnarray*}
a(\phi^*(r_2))-a(\phi^*(r_1)) &\leq & M_1 \mu (\{\phi^*(r_2)\leq \phi^*(|y|)\leq \phi^*(r_1)\})\\
&=&M_1 \mu (r_1 \leq |y|\leq r_2)=M_1 \pi (r_2^2-r_1^2)\\
&\leq & 2\pi R M_1  (r_2-r_1). 
\end{eqnarray*}
On the other hand, 
\begin{eqnarray*}
a(\phi^*(r_2))-a(\phi^*(r_1)) &\geq & m \mu (\{\phi^*(r_2)\leq \phi(y)\leq \phi^*(r_1)\})\\
&\geq& \frac{m}{M_2} \int_{\{\phi^*(r_2)\leq \phi(y)\leq \phi^*(r_1)\}} |\nabla \phi| dy \\
&\geq & \frac{m}{M_2} \int_{\phi^*(r_2)}^{\phi^*(r_1)} \int_{\{\phi^{-1}(t)\}}dsdt \\
&\geq & \frac{m}{M_2}(\phi^*(r_1)-\phi^*(r_2)) K(r_1,r_2),
\end{eqnarray*}
where 
\[K(r_1,r_2)=\min \limits_{\{ \phi^*(r_2) \leq t \leq \phi^*(r_1)\}} \mathcal{H}^{n-1}(\phi^{-1}(t))>0, \ \ 0<r_1< r_2 <  R.\] 
Since 
$\{\phi^{-1}(t)\}= \partial \{x: \phi(x)>t\}$, it follows from the isoperimetric inequality that if $\phi^*(r_1)< \max \limits_{y \in \overline{\Omega}} \phi -\delta$ and $\phi^*(r_2)> \min \limits_{y \in \overline{\Omega}} \phi +\delta$, for some $\alpha \in (0,1)$, then  \[K(r_1,r_2)>C>0, \ \ \forall r_2 \ \ \hbox{with}\ \  r_1<r_2 < R,\]
for some $C>0$ independent of $\phi$. Hence we have
\begin{eqnarray}
0\leq \frac{a(\phi^*(r_2))-a(\phi^*(r_1))}{r_2-r_1}\leq \frac{2\pi R M_1M_2}{mK(r_1,r_2)} \leq \frac{2\pi R M_1M_2}{mC}. 
\end{eqnarray}
By approximation the above also holds for $\epsilon<r_1<r_2< R-\epsilon$. Thus $\phi^*$ is Lipschitz continuous on $(\epsilon, R-\epsilon)$ for every $\epsilon>0$.

\hfill $\Box$

Now let
\[
j(t):=\int_{\{\phi >t \} } |\nabla \phi |^2 dy, \quad j^*(t):=\int_{\{\phi >t \} } |\nabla \phi^* |^2 dy,  \ \ \forall  t>\min_{y\in \overline \Omega}\phi;
\]
\[J(t):=\int_{\{\phi>t\}}|\nabla \phi| dy, \quad J^*(t):=\int_{\{\phi^*>t\}}|\nabla \phi^*| dy, \ \ \forall  t>\min_{y\in \overline \Omega}\phi.
\]
It is easy to see that both $j(t)$ and $J(t)$ are  absolutely continuous  and decreasing  in $ t>\min_{y\in \overline \Omega}\phi $.  If $ \phi \equiv C $ on $\partial \Omega$,  it
can  be shown that 
\begin{equation}\label{gradient}
\int_{\{\phi=t\}}|\nabla \phi| ds\geq \int_{\{\phi^*=t\}}|\nabla \phi^*|ds, \quad  \hbox{for  \it {a.e.}}  \ \  t>\min_{y\in \overline \Omega}\phi.
\end{equation}
Indeed it follows from  Cauchy-Schwarz and Bol's inequalities that 

 \begin{eqnarray*}
 \int_{\{\phi=t\}}|\nabla \phi| ds
&\geq & \left(\int_{\{\phi=t\}} e^{\frac{u}{2}} \right)^{2} \left( \int_{\{\phi=t\}} \frac{e^{u}}{|\nabla \phi|}\right)^{-1}\\
&=& \left(\int_{\{\phi=t\}} e^{\frac{u}{2}} \right)^{2} \left( -\frac{d}{dt}\int_{\Omega_t}e^{u}\right)^{-1}\\
&\geq &\frac{1}{2}(\int_{\Omega_t}e^{u})(8\pi - \int_{\Omega_t}e^{u})(-\frac{d}{dt}\int_{\Omega_t}e^{u})^{-1}\\
&=&  \frac{1}{2}(\int_{\Omega^*_t}e^{U_{\lambda}})(8\pi - \int_{\Omega^*_t}e^{U_{\lambda}})(-\frac{d}{dt}\int_{\Omega^*_t}e^{U_{\lambda}})^{-1}\\ 
 &=& \int_{\{\phi^*=t\}}|\nabla \phi^*| ds, \ \ \hbox{for  \it {a.e.}}  \ \ t>\min_{y\in \overline \Omega}\phi.
 \end{eqnarray*}
It also follows that   $j^*(t),  J^*(t)$ are  absolutely  continuous and decreasing  in $t>\min_{y\in \overline \Omega}\phi$,  since both  functions are right-continuous by definition and 
 \[
 0\le j^*(t-0)-j^*(t)\le j(t-0) -j(t)=\int_{\{\phi=t\}}|\nabla \phi|^2 dy =0,   \ \ t>\min_{y\in \overline \Omega}\phi. \]
 \[
 0\le J^*(t-0)-J^*(t)\le J(t-0) -J(t)=\int_{\{\phi=t\}}|\nabla \phi|dy =0,   \ \ t>\min_{y\in \overline \Omega}\phi. \]
 Therefore we have  the following proposition. 
 
\begin{proposition}\label{rearrangProp}
Let $u\in C^2(\overline{\Omega})$ satisfy 
\[\Delta u+e^{u}\geq 0\ \ \hbox{in}\ \ \Omega,\]
and let $U_{\lambda}$ be given by (\ref{ULambda}). Suppose $\phi \in C^1(\overline{\Omega})$ and $\phi \equiv C$ on $\partial \Omega$. 
Define the equimeasurable symmetric rearrangement $\phi^*$ of $\phi$, with respect to the measures $e^{u}dy$ and $e^{U_\lambda}dy$, by (\ref{rearrang}). Then  $\phi^*$ is Lipschitz continuous on $(\epsilon,R-\epsilon)$ for every $\epsilon>0$, and   $j^*(t), J^*(t)  $ are  absolutely continuous and decreasing  in $t>\min_{y\in \overline \Omega}\phi$  and \eqref{gradient} holds.\\ 
\end{proposition}

\section{Bol's type  inequalities}

We first prove the following lemma.

\begin{lemma}\label{goodLemma!}
Let $\psi \in C(\R^2\setminus B_R)$ be a decreasing radial function and 
\[\int_{(\R^2\setminus B_R)}e^{\psi}dx <\infty,\]
for some $R>0$. Then 
\[\lim_{s\rightarrow -\infty} e^{s}\int_{\{\psi >s\}}dx=0.\]
\end{lemma}
{\bf Proof.} Since $\psi$ is decreasing, 
\begin{eqnarray*}
\frac{3 \pi}{4} r^2e^{\psi (r)} \leq \int_{(B_r \setminus B_{r/2})}e^{\psi}dx, 
\end{eqnarray*}
for $r>2R$. Letting $r\rightarrow \infty$ we obtain, 
\[\lim_{r \rightarrow \infty}r^2e^{\psi (r)} =0.\]
Define
\[ r(s):= \sup\{ r \ge R:  \psi(r) >s\}, \quad s \in \R.  \]
Then $r(s) $  is well-defined for $ s <  \psi(R)$ and $ \lim_{ s \to -\infty} r(s) =\infty$.
Since
\[ e^{s}\int_{\{\psi >s\}}dx \le    \pi ( r(s)^2- R^2) e^{\psi( r(s) )},   \]
we obtain
\[\lim_{s\rightarrow -\infty} e^{s}\int_{\{\psi >s\}}dx=0.\]
The proof is complete. \hfill $\Box$ \\ \\
For the proof of our main results, we shall need the following reversed Bol's inequality. 

\begin{proposition}\label{RadialBolPropExterior}
Let $B_R$ be the ball of radius $R$ in $\R^2$ $\psi \in C^{0,1}(\R^2 \setminus B_R)$ be a strictly decreasing radial function satisfying
\begin{equation}\label{MainPDECondRelaxedExterior}
\int_{\partial B_r}|\nabla \psi |ds \leq 8\pi - \int_{\R^2 \setminus B_r}e^{\psi}  \ \ \hbox{for a.e. } r\in(R,\infty), \ \ \hbox{and} \ \ \int_{\R^2 \setminus B_R}e^{\psi}< 8 \pi.
\end{equation}
Then the following inequality holds
\begin{equation}\label{RadialBolInequalityExterior}
\left( \int_{\partial B_R}e^{\frac{\psi}{2}}\right)^2\leq \frac{1}{2}\left( \int_{\R^2 \setminus B_R} e^{\psi}\right) \left( 8 \pi-\int_{\R^2 \setminus B_R}e^{\psi}\right).
\end{equation}
Moreover if $\int_{\partial B_r}|\nabla \psi|ds \not \equiv  8\pi - \int_{\R^2 \setminus B_r}e^{\psi} $ on $(R, \infty)$, then the inequality in (\ref{RadialBolInequalityExterior}) is strict. 
\end{proposition}
{\bf Proof.} Let $\beta:=\psi(R)$ and define
\[k(s)=8\pi - \int_{\{\psi<s\}}e^{\psi}dx, \ \ \hbox{and} \ \ \mu(s)=\int_{\{\psi>s\} }dx+\pi R^2,\]
for $s<\beta$. Then 
\[-k'(s)=\int_{\{\psi=s\}}\frac{e^{\psi}}{|\nabla \psi|}=-e^{s}\mu'(s).\]
Hence 
\begin{eqnarray}\label{RelaxedBolSharp}
-k(s)k'(s)&\geq& \int_{\{\psi=s\}}|\nabla \psi | \cdot \int_{\{\psi=s\}}\frac{e^{\psi}}{|\nabla \psi|}\\
&= & (\int_{\{\psi=s\}}e^{\psi/2})^{2}=e^{s}(\int_{\{\psi=s\}}ds)^{2}\nonumber \\
&= & e^s  \cdot 4 \pi (\int_{\{\psi>s\}}dx+\pi R^2)=4 \pi e^s\mu(s),\nonumber
\end{eqnarray}
for a.e. $s< \beta$.
Therefore 
\[\frac{d}{ds} [e^s \mu(s)-k(s)+\frac{1}{8\pi}k^2(s)]=\mu(s)+\frac{1}{4\pi}k'(s)k(s)\leq 0,\]
for a.e. $s<\beta$. Integrating on $(- \infty, \beta)$ and using Lemma \ref{goodLemma!} we get 
\begin{equation}\label{calInequality}
\left[ e^s \mu(s)-k(s)+\frac{1}{8\pi}k^2(s)\right]_{-\infty}^{\beta}=e^{\beta} \mu(\beta)-k(\beta)+\frac{1}{8\pi}k^2(\beta)\leq 0.
\end{equation}
Now notice that 
\[ k(\beta)=8 \pi- \int_{\R^2\setminus B_R}e^{\psi}dx\]
and
\[e^{\beta}\mu(\beta)=e^{\beta} \int_{ B_R}dx = \frac{1}{4\pi} e^{\beta}(\int_{\partial B_R}ds)^{2}=\frac{1}{4 \pi}(\int_{\partial B_R} e^{\frac{\psi}{2}}ds)^{2}. \]
Thus (\ref{RadialBolInequalityExterior}) follows from the inequality (\ref{calInequality}). Finally if  $\int_{\partial B_r}|\nabla \psi|ds \not \equiv 8\pi - \int_{\R^2 \setminus B_r}e^{\psi}  $ on $(R, \infty)$, then the inequality (\ref{RelaxedBolSharp}) will be strict, and consequently (\ref{RadialBolInequalityExterior}) will also be strict. \hfill $\Box$\\ \\
Similarly one can prove the following proposition (see Proposition 2.2 in \cite{G-M1}).

\begin{proposition}\label{RadialBolProp}
Let $B_R$ be the ball of radius $R$ in $\R^2$ and  $u\in C^{0,1}(\overline{B_R})$ be a strictly decreasing radial function satisfying
\begin{equation}\label{MainPDECondRelaxed}
\int_{\partial B_r}|\nabla u|ds \leq \int_{B_r}e^{u}   \ \ \hbox{for a.e. } \ \ r\in(0,R), \ \ \hbox{and} \ \ \int_{B_R}e^{u}< 8 \pi.
\end{equation}
Then the following inequality holds
\begin{equation}\label{RadialBolInequality}
\left( \int_{\partial B_R}e^{\frac{u}{2}}\right)^2\geq \frac{1}{2}\left( \int_{B_R} e^u\right) \left( 8 \pi-\int_{B_R}e^u\right).
\end{equation}
Moreover if $\int_{\partial B_r}|\nabla \psi|ds \not \equiv \int_{B_r}e^{\psi} $ on $(0, R)$, then the inequality in (\ref{RadialBolInequality}) is strict. 
\end{proposition}

\section{Proof of the main results} 

\begin{lemma}\label{LastEstimateExterior}
Let $R>0$ and assume that $\psi \in C^{0,1}(\R^2 \setminus B_R)$ is a strictly decreasing radial function that satisfies 
\begin{equation}\label{superSol1}
\int_{\partial B_r} |\nabla \psi|  \leq 8\pi -\int_{\R^2 \setminus B_r}e^{\psi}
\end{equation}
for a.e. $r\in (R,\infty)$ and $\psi=U_{\lambda_1}=U_{\lambda_2}$ on $\partial B_R$ for some  $\lambda_2> \lambda_1$. Then 
\begin{equation}\label{lastEstimate1}
\int_{\R^2 \setminus B_R} e^{U_{\lambda_2}} \leq  \int_{\R^2 \setminus B_R} e^{\psi} \leq  \int_{\R^2 \setminus B_R} e^{U_{\lambda_1}}.
\end{equation}
Moreover if $\int_{\partial B_r} |\nabla \psi|  \not \equiv \int_{\R^2 \setminus B_r}e^{\psi}$ on $r\in(R, \infty)$, then the inequalities in (\ref{lastEstimate1}) are also strict.
\end{lemma}
{\bf Proof.} Let $m_1:=\int_{\R^2 \setminus B_R}e^{U_{\lambda_1}}$, $m_2:=\int_{\R^2 \setminus B_R}e^{U_{\lambda_2}}$, and $m:=\int_{\R^2 \setminus B_R}e^{\psi}$. Define 
\[\beta:=\left( \int_{\partial B_R} e^{\frac{\psi}{2}}\right)^2=\left( \int_{\partial B_R} e^{\frac{U_{\lambda_1}}{2}}\right)^2=\left( \int_{\partial B_R} e^{\frac{U_{\lambda_2}}{2}}\right)^2\]
It follows from Proposition \ref{RadialBolPropExterior} that 
\[\beta\leq \frac{1}{2}m(8 \pi -m)\]
On the other hand 
\[\beta=\frac{1}{2}m_1(8 \pi -m_1)=\frac{1}{2}m_2(8 \pi -m_2),\]
i.e. $m_1$ and $m_2$ are roots of the quadratic equation
\[x^2-8\pi x +2 \beta=0.\]
Since $m$ satisfies
\[m^2-8\pi m +2\beta \leq 0,\]
we have 
\[m_2 \leq m \leq m_1.\]
\hfill $\Box$ \\ \\ 
Similarly the following lemma holds (see Lemma 3.3 in \cite{G-M1}). 

\begin{lemma}\label{LastEstimate}
Assume that $\psi \in C^{0, 1} (\overline{B_R})$ is a strictly decreasing, radial, Lipschitz function,  and  satisfies 
\begin{equation}\label{superSol}
\int_{\partial B_r} |\nabla \psi|  \le \int_{B_r}e^{\psi}
\end{equation}
{\it a.e.}  $r\in (0,R)$ and $\psi=U_{\lambda_1}=U_{\lambda_2}$  for some $\lambda_2> \lambda_1$ on $\partial B_R$, and $R>0$. Then 
there holds 
\begin{equation}\label{lastEstimate}
\hbox{either }  \,\,  \int_{B_R} e^{\psi} \le  \int_{B_R} e^{U_{\lambda_1}}  \,\,\quad  \hbox{or} \,\, \quad  \int_{B_R} e^{\psi}\geq   \int_{B_R} e^{U_{\lambda_2}}.
\end{equation}
Moreover if the inequality in (\ref{superSol}) is strict in a set with positive measure in $(0,R)$, then the inequalities in (\ref{lastEstimate}) are also strict.
\end{lemma}

We shall also need the following lemma. \\
\begin{lemma}\label{RadialComparision}
Assume that $\psi \in C^{0,1}(\overline{B_R})$ is a strictly decreasing radial function satisfying \eqref{superSol} for a.e. $r\in (0,R)$. If 
\[\rho=\int_{B_R}e^{\psi}dx =\int_{B_R}e^{U_\lambda}< 8\pi ,\]
then $U_\lambda (R) \leq \psi(R) $.

\end{lemma}
{\bf Proof.} By Proposition \ref{RadialBolProp} we have 
\begin{eqnarray*}
\left( \int_{\partial B_R}e^{\frac{U_\lambda}{2}}\right)^2 &=& \frac{1}{2}\left( \int_{B_R} e^{U_\lambda}\right) \left( 8 \pi-\int_{B_R}e^{U_\lambda}\right)\\
&= & \frac{1}{2}\left( \int_{B_R} e^{\psi}\right) \left( 8 \pi-\int_{B_R}e^{\psi}\right)\\
&\leq& \left( \int_{\partial B_R}e^{\frac{\psi}{2}}\right)^2, 
\end{eqnarray*}
and hence $U_\lambda (R) \leq \psi (R)$. \hfill $\Box$ \\ \\
Now we are ready to prove the main result of this paper, Theorem \ref{New-8Pi-LowerBound}. \\ \\
{\bf Proof of Theorem \ref{New-8Pi-LowerBound}.} First we prove that $\rho \geq 8\pi$. Suppose $w_1$ and $w_2$ satisfy the assumptions of Theorem \ref{New-8Pi-LowerBound}. Then 
\[\Delta (w_2-w_1)+e^{w_2}-e^{w_1}=f_2-f_1\geq 0.\]
Now  choose $\lambda>0$ and $R\in (0,\infty)$ such that 
\begin{equation}
\int_{\Omega}e^{w_1}=\int_{B_R}e^{U_{\lambda}},
\end{equation}
and let $\phi$ be the symmetrization of $w_2-w_1$ with respect to the measures $e^{w_1}dy$ and $e^{U_{\lambda}}dy$. Then it follows from Proposition \ref{rearrangProp} and Fubini's theorem that
\begin{eqnarray*}
\int_{\{\phi=t\}}|\nabla \phi| &\leq &  \int_{\{w_2-w_1=t\}}|\nabla (w_2-w_1)| \\
&\leq&\int_{\Omega_t}e^{w_2}-e^{w_1}dx\\
&=& \int_{\{\phi>t\}}e^{U_{\lambda}+\phi}- \int_{\{\phi>t\}}e^{U_{\lambda}}\\
&=&\int_{\{\phi>t\}}e^{U_{\lambda}+\phi}- \int_{\{\phi=t\}} |\nabla U_{\lambda}|,
\end{eqnarray*} 
 for a.e. $t>\inf_{\Omega}(w_2-w_1)$. Hence
\begin{equation}
\int_{\{\phi=t\}} |\nabla (U_{\lambda}+\phi)| \leq \int_{\{\phi>t\}} e^{(U_{\lambda}+\phi)} d
\end{equation}
for all $t>\inf_{\Omega}(w_2-w_1)$.  Since  $\phi  $ is decreasing in $r$,     $\psi:= U_{\lambda}+\phi$ is a strictly decreasing function, and 
\begin{equation}\label{supersolution}
\int_{\partial B_r} |\nabla \psi| \le  \int_{B_r} e^{\psi} dy, \quad {\it a.e.}  \quad r \in (0, R), 
\end{equation}
by Proposition \ref{rearrangProp} and the above inequality we see that  $\psi \in W^{1,\infty}(B_R)$ and thus by Morrey's inequality $\psi \in C^{0, 1} (B_R)$.

Since $w_1 \not \equiv w_2$ and $\int_{\Omega}e^{w_1}=\int_{\Omega}e^{w_2}$, then $w_2<w_1$ on a subset of $\Omega$ with positive measure. Hence  $\phi(R)<0$ and consequently $\psi(R)=U_{\lambda}(R)+\phi(R)<U_{\lambda}(R)$. This is a contradiction in view of Lemma \ref{RadialComparision}, and therefore we must have $\rho\geq 8\pi$.  

Next we prove that $\rho >8 \pi$. Suppose $\rho=8\pi$ and let $\lambda_1>0$. With an argument similar to the one above we may show that there exists $\psi=U_{\lambda_1}+\phi \in C^{0,1}(\R^2)$ such that 
\[\int_{\Omega}e^{w_1}dx=\int_{\R^2}e^{U_{\lambda_1}}=8\pi= \int_{\Omega}e^{w_2}dx=\int_{\R^2}e^{\psi}dx,\]
and 
\begin{equation}\label{psiDiffInequality}
\int_{\partial B_r}|\nabla \psi| \leq \int_{B_r}e^{\psi}dx
\end{equation}
for a.e. $r \in (0,\infty)$. Since $\int_{\R^2}e^{\psi}=\int_{\R^2}e^{U_{\lambda_1}}$, there exists $r_0\in(0,\infty)$ such that $\psi(r_0)=U_{\lambda_1}(r_0)$.  There exists a positive constant $\lambda_2\neq \lambda_1$ such that $U_{\lambda_2}(r_0)=U_{\lambda_1}(r_0)=\psi(r_0)$. Since $\psi > U_{\lambda_1}$ in $B_{r_0}$, it follows from Proposition \ref{LastEstimate} that $\lambda_1<\lambda_2$ and 
\[\int_{B_{r_0}}e^{\psi} \geq \int_{B_{r_0}}e^{U_{\lambda_2}}.\] 
On the other hand $\psi< U_{\lambda_1}$ in $\R^2 \setminus B_{r_0}$ and consequently it follows from Proposition $\ref{LastEstimateExterior}$ that 
\[\int_{\R^2 \setminus B_{r_0}}e^{\psi}\geq \int_{\R^2 \setminus B_{r_0}}e^{U_{\lambda_2}}.\]
Hence 
\begin{equation}\label{comparision0}
8 \pi =\rho = \int_{\R^2}e^{\psi}\geq \int_{\R^2 }e^{U_{\lambda_2}}=8\pi.
\end{equation}
 Note that if $f_1 \not \equiv 0$ or $f_2  \not \equiv f_1$, then  the inequality in \eqref{comparision0} will be strict,  which is a contradiciton . Suppose $f_1\equiv f_2\equiv 0$.  We may assume without loss of generality that $c= w_2-w_1\ge 0$ on $\partial \Omega$, since otherwise we can switch $w_1$ and $w_2$.     By \eqref{comparision0}  we conclude that   the equality  in  \eqref{psiDiffInequality} holds for  a.e. $r \in (0,\infty)$ and $\psi =U_{\lambda_2}$.  It also  yields that  the equality in  \eqref{gradient}  must be true for $\phi=w_2-w_1$  and  $t \ge \inf_{\Omega} \phi$.  By  the proof of Proposition \ref{rearrangProp},   we  also know that  Bol's inequality  \eqref{BolInequality0} on $\omega =\{ \phi >t \}$  must be equality,   and therefore $\{\phi >t \}$ must be simply-connected  for $t \ge \inf_{\Omega} \phi$  by Theorem D.   This is  a contradiction
 since $\{\phi >t \}$ is not  simply-connected when $  \inf_{\Omega} \phi <t<0 $.   The contradiction implies  $\rho> 8 \pi$.
  \hfill $\Box$\\
\hspace{.5cm}

{\bf Proof of Theorem \ref{New-8Pi-LowerBoundM}.} The proof follows from Theorem D and the same argument used in the proof of Theorem \ref{New-8Pi-LowerBound}.  \hfill $\Box$ \\

$\mathbf{Acknowledgement}$  The authors would like to thank the anonymous referee for their careful reading of the paper and many helpful comments. The first author is partially supported by a Simons Foundation Collaborative Grant (Award \#199305) and NSFC grant No 11371128. The second author is partially supported by a start-up grant from University of California, Riverside.


\begin{thebibliography}{12}

\bibitem{Aubin} T. Aubin, Mulleures constantes dans des theoremes d'inclusion de Sobolev at un
theoreme de Fredholm non-linearire pour la transformation conforme de la courbure
scalaire. J. Funct. Anal., 32 (1979), 149-179.

\bibitem {Bandle} C. Bandle, Isoperimetric Inequalities and Applications. London: Pitman (1980).


\bibitem{BLT}  D. Bartolucci, C.S.  Lin, G.  Tarantello, Uniqueness and symmetry results for solutions of a mean field equation on $S^2$ via a new bubbling phenomenon. Comm. Pure Appl. Math. 64 (2011), no. 12, 1677-1730. 

\bibitem{BL0}  D. Bartolucci, C.S.  Lin,   Uniqueness results for mean field equations with singular data. Comm. Partial Differential Equations 34 (2009), no. 7-9, 676-702. 

\bibitem{BL}  D. Bartolucci, C.S.  Lin,  Existence and uniqueness for mean field equations
on multiply connected domains at the critical parameter. Math. Ann. 359 (2014), 1-44.

\bibitem{Bol} G. Bol, Isoperimetrische Ungleichungen fur Bereiche auf Fl\"achen, Jahresber. Deutschen Math. Vereinigung 51 (1941), 219-257. 

\bibitem{BZ} Y.D. Burago, V.A. Zalgaller, Geometric Inequalities, Springer Ser. Sov. Math., Springer-Verlog (1988).

\bibitem{CCL} S.Y.A. Chang, C.C. Chen, C.S. Lin, Extremal functions for a mean field equation in two dimension.  Lecture on Partial Differential Equations, New Stud. Adv. Math., vol. 2, pp. 61–93. Int. Press, Somerville (2003).

\bibitem{CY1}  S.Y.A. Chang, P. Yang, Conformal deformation of metrics on $S^2$.   J. Diff. Geom. 27(2), 259-296 (1988). 

\bibitem{CY2}  S.Y.A. Chang, P. Yang, Prescribing Gaussian curvature on $S^2$,  Acta Math. 159(3-4), 215-259 (1987).

\bibitem{ChengLin}  K.S. Cheng, C.S. Lin, On the asymptotic behavior of solutions of the conformal Gaussian curvature equations in $\R^2.$  Math. Ann. 308(1), 
119-139 (1997). 

\bibitem{FFGG} Feldman, J., Froese, R., Ghoussoub, N., Gui, C.F.: An improved Moser-Aubin-Onofri inequality for axially symmetric functions on S2. Calc. Var. Part. Diff. Eqs. 6(2), 95-104 (1998).

\bibitem {GL} N. Ghoussoub; C.S. Lin, On the best constant in the Moser-Onofri-Aubin inequality. Comm. Math. Phys. 298 (2010), no. 3, 869-878. 

\bibitem{G-M1} C. Gui and A. Moradifam,  The Sphere Covering Inequality and its applications, preprint.

\bibitem{GW}  C. Gui, J.C. Wei, On a sharp Moser-Aubin-Onofri inequality for functions on $S^2$ with symmetry.  Pac. J. Math. 194(2), 349-358 (2000) 


\bibitem {GMbook} N. Ghoussoub, A. Moradifam, Functional Inequalities: New Perspectives and Applications. Mathematical Surveys and Monographs, American Mathematical Society, 2013.  

\bibitem {Matano} H. Matano, Nonincrease of the lap-number of a solution for a one-dimensional semilinear parabolic equation. J. Fac. Sci. Univ. Tokyo Sect. IA Math. 29 (1982), no. 2, 401-441. 

\bibitem{Moser} J. Moser, A sharp form of an Inequality by N. Trudinger. Ind. Univ. Math. J.  20 (1971), 1077-1091.


\bibitem{Lin} C.S. Lin, Uniqueness of solutions to the mean field equations for the spherical Onsager vortex. Arch. Rat. Mech. Anal. 153(2), 153-176 (2000).

\bibitem{LinDuke} C.S. Lin, Topological degree for mean field equations on $S^2$. Duke Math J. 104(3) 501-536 (2000),.

\bibitem {LL} C.S. Lin, M. Lucia, One-dimensional symmetry of periodic minimizers for a mean field equation. Ann. Sc. Norm. Super. Pisa Cl. Sci. (5) 6 (2007), no. 2, 269-290.

\bibitem{Onofri} E. Onofri, On the positivity of the effective action in a theorem on radom surfaces.
Comm. Math. Phys., 86 (1982), 321-326.


\bibitem{Suzuki} T. Suzuki, Global analysis for a two dimensional elliptic eigenvalue problem with exponential nonlinearity. Ann. I.H.P. Analyse Non Lin\'{e}aire 9 (1992), 367-398. 

\end{thebibliography}
 \end{document}